\documentclass[11pt]{article}

\usepackage{geometry}                
\usepackage{graphicx}
\usepackage{amssymb}
\usepackage{epstopdf}

\title{Shifted small deviations and Chung LIL for symmetric alpha-stable processes.}
\author{Elena Shmileva\thanks{Financial support from the Austrian Science Fund (FWF) under the grant P18022 and START-project Y328 is gratefully acknowledged.}}
\date{}                                           

\def\E{\mathbf{E}\,}
\def\Pr{\mathbf{P}\,}

\newtheorem{theorem}{Theorem}
\newtheorem{lemma}{Lemma}
\newtheorem{prop}{Proposition}

\newtheorem{fact}{Fact}
\newtheorem{corol}{Corollary }

\newenvironment{proof}{{\it{Proof.}\qquad}}{ \begin{flushright} $\blacksquare$ \end{flushright}}

\begin{document}

\maketitle

\begin{abstract}
Let $X_{\alpha}$ be a symmetric $\alpha$-stable L\'evy process with $\alpha\in (1,2)$. We consider small ball probabilities of the following type $\Pr\{ \|X_{\alpha}-\lambda\, f\|< r \}$ as $r\to 0$ and $\lambda r^{\alpha-1}\to 0$ or $\lambda r^{\alpha-1} = c$, $c>0$, where $\|\cdot\|$ is the sup-norm and $f$ is any continuous function which starts at $0$. We obtain an exact rate of decrease for these probabilities including constants.

Using these small ball estimates, we derive a functional LIL for $X_{\alpha}$ with continuous attracting   functions. It occurs that the a.s. limit set of the family $\left\{\frac{X_{\alpha}(T\cdot)}{T^{1/\alpha}h(T)}\right\}_{T>0}$  is equal to the set of all continuous functions (which start at $0$), under certain choice of scaling function $h(T)$. 

\small{Keywords:\quad Chung LIL, Strassen LIL, small ball probabilities, L\'evy processes, stable L\'evy processes.} 
\end{abstract}

\section*{Introduction}

We are interested in the probabilities that a c\`adl\`ag process $X(t),
t\in[0,1]$ hits an arbitrary small ball, i.e., $\Pr\{X\in
\mathrm{B}(f,r) \}$, where  $\mathrm{B}(f,r)$ is  a ball (in the
Skorokhod metric or in the uniform metric) of radius $r>0$, ($r\to 0$), and of center $f$, which  is an arbitrary element of the Skorokhod space  $D[0,1]$.

If the shift function (center) $f$ has jumps, i.e., $f\in D[0,1]
\setminus C[0,1] $, then the problem is delicate. If the process $X$ has no fixed-time jumps, what holds in most of practically important cases, the uniform small
balls are obviously empty. Hence, one has
to deal with the Skorokhod topology. We don't know any results on
probabilities of small balls in the Skorokhod topology. This is a
subject of future research.

However, if we assume that  the shift  $f$ is a continuous
function, then there is a sense to consider small balls in the
uniform topology (as well as in the Skorokhod topology).  Dealing
with uniform balls  and the uniform topology is more usual and
there are already some results in this direction.

In the sequel,  by $\|\cdot\|$ we denote the uniform norm, and by $B(f, r)$ a ball of radius $r$ and of center $f$ in the uniform metric.

Aurzada and Dereich (see  \cite{AD_08})  elaborate a method that allows to estimate  $\Pr\{\|X\|< r\}$, where $X$ is an arbitrary L\'evy  process. 
So, they deal with  time-homogeneous processes, i.e., the shift  $f$ is the identity function multiplied by a constant.  Since we are interested in applications to the functional law of the iterated logarithm (functional LIL), we need to study similar probabilities but with arbitrary shift functions. Thus, in general, we deal with time-inhomogeneous (additive) processes.

In this article, we focus on symmetric $\alpha$-stable L\'evy processes $X_{\alpha}$. 

Concerning the LIL for these processes, there are significant differences from
the gaussian and the pre-gaussian cases. 
Namely, Limsup LIL doesn't exist, i.e., there is no such a scaling function $\varphi(\cdot)$ that
$0<\limsup_{t\to\infty}{|X_{\alpha}(T)|}/{\varphi(T)}<\infty$. Instead, there is an integral test for  $\varphi$ (see Fact \ref{fact 2} below) that says whether this limit is equal to $0$ or to $\infty$.  

In spite of that, there is  a Liminf LIL statement by Taylor \cite{Tay_67}:
\begin{eqnarray*}
\liminf_{T\to \infty}\frac{\|X_{\alpha}(T\cdot)\| }{(T/\log\log T)^{1/\alpha} }=K_{\alpha}^{1/\alpha} \textit{\ \ a.s.,\ \ }
\end{eqnarray*}
where $K_{\alpha}$ is a positive constant (the same as in (\ref{Mog}) below).
 
Based on these two facts, we are looking for a functional LIL  for $X_{\alpha}$ under those scaling  functions $\varphi$, which are bigger than  $(T/\log\log T)^{1/\alpha} $.  
For example, if $\ \varphi(T) \cdot (T/\log\log T)^{-1/\alpha}\to \infty$, then  
\begin{eqnarray*}
\liminf_{T\to \infty}\left\|\frac{X_{\alpha}(T\cdot) }{\varphi(T)}\right\|=0 \textit{\ \ a.s.,\ \ }
\end{eqnarray*}
what means that the family of scalings $\left\{\frac{X_{\alpha}(T\cdot)}{\varphi(T)}\right\}_{T>0}$ has at least one a.s. limit point under uniform convergence, this is the zero function.
If, moreover, the integral test gives $0$, then this is the only a.s. limit point.  

In this article, we study the a.s. limit sets of the family under these scaling functions $\varphi$ that ensure $\infty$ in the integral test.    
In Theorems \ref{LILL} and  \ref{LIL},  we obtain that if $\varphi(T)\in ((T/\log\log T)^{1/\alpha}, T^{1/\alpha}\log\log T^{1-1/\alpha})$, then the a.s. limit set of $\left\{\frac{X_{\alpha}(T\cdot)}{\varphi(T)}\right\}_{T>0}$ in the uniform topology is equal to the set of all  continuous functions which start at $0$.  

The border line $\varphi(T)= C\cdot (T/\log\log T)^{1/\alpha}, C>0$ is studied in Theorem \ref{Chung_deg}, which shows that the scaling is too small and the trajectories stop a.s. clustering near continuous functions, i.e., the a.s. limit set is empty. 

Of course, it is interesting to understand what happens when the scaling function is close to the border of the integral test. It requires additional study. 

The article is structured as follows.

In section 1, we obtain small deviation estimates for  $\Pr\{ \|X_{\alpha}-\lambda\, f\|< r \}$ as $r\to 0$, first under $\lambda r^{\alpha-1}\to 0$, see Theorem \ref{small lambda}, then under $\lambda r^{\alpha-1} = c$ , $c>0$, see Theorem \ref{LIL case}. For the proof we use the Girsanov theorem for additive processes, that is provided in section 0.1.

We start section 2 with a detailed review on the LIL for stable L\'evy processes and  discuss  a Baldi-Royonette result (see \cite{BR_92}) for the Wiener process that describes a parallel situation with the main results of this article.
In Theorems \ref{Chung_deg},  \ref{LILL} and  \ref{LIL}, we get the a.s. limit set (subset) of the families $\left\{\frac{X_{\alpha}(T\cdot)}{T^{1/\alpha}h(T)}\right\}_{T>0}$, when $h(T)\in[(\log\log T)^{-1/\alpha}, (\log\log T)^{1-1/\alpha}]$. We obtain these results together with the rates of convergence to the limit functions.

\subsection{Notations and  tools.} 

Let ${\cal{C}}= \{f\in C[0,1] : f(0)=0\}$.
\\
By $AC[0,1]$ we denote the set of all absolutely continuous functions on $[0,1]$.
\\
Put ${\bf{H}}= \{f\in AC[0,1] :\ f'\in L_2, \  f(0)=0\}$.
\\
 We use notations from the Sato monograph
\cite{Sat_99}, to introduce a L\'evy process. In this article,  we  deal with  processes of finite
expectation, therefore  it is convenient to define  a L\'evy process $X$
by its centered triplet $(\sigma^2, \Lambda(dx), \gamma )_1$,
where $\sigma^2$ is the variance of  the gaussian component (here, $\sigma^2=0$), $\Lambda(dx)$ is the L\'evy measure and $\gamma$
is the expectation of $X(1)$. Here is the corresponding L\'evy-Ito decomposition
\begin{eqnarray*}
X(t)=\gamma t+\int_0^t \int_{{\mathbb{R}}\setminus \{0\}}x\bar{\cal{N}}_{\Lambda}(dx,ds),
\end{eqnarray*}
were ${\cal{N}}_{\Lambda}(dx, dt)$ is a Poisson measure corresponding to the L\'evy measure $\Lambda$, and $\bar{\cal{N}}_{\Lambda}(dx, dt)={\cal{N}}_{\Lambda}(dx, dt)-\Lambda(dx) dt$ is the centered Poisson measure. 

If $\gamma=0$, we call the corresponding L\'evy process  a $(\Lambda, 0)$-L\'evy martingale.

By additive processes, we mean time-inhomogeneous  processes with independent increments, that start at $0$. The distributions  of the processes with finite expectations are specified by the centered triplets $(0, \Lambda(dx, dt), \gamma(t) )_1$, $\gamma\in L_1$. The corresponding L\'evy-Ito decomposition is
\begin{eqnarray*}
X(t)=\int_0^t \gamma(s)ds+\int_0^t \int_{{\mathbb{R}}\setminus \{0\}}x\bar{\cal{N}}_{\Lambda}(dx,ds).
\end{eqnarray*}
Denote by $P_{\xi}$ the distribution  of the process $\xi$ in $D[0,1]$.

In the next section, we will need the following particular case of the Girsanov theorem, see Theorem 3.24 from \cite{JS_03}, see also \cite{LS_02}, Theorem 2:
\begin{fact}$($The Girsanov transform for additive processes with finite expectations$)$

Let $\xi$  be an additive process defined by the centered triplets
$(0, \Lambda(dx, dt), \gamma(t))_1$, $\gamma\in L_1[0,1]$. 
Suppose there exists $\theta(\cdot, \cdot): {\mathbb{R}}\setminus \{0\}\times [0,1]\to \mathbb{R}$ such that
\begin{eqnarray}\label{cond}
\int_0^1 \int_{{\mathbb{R}}\setminus \{0\}} \left(e^{\theta(x,s)/2}-1\right)^2\Lambda(dx,ds)<\infty.
\end{eqnarray}
Then the distribution of an additive process $\xi_{\theta}$ defined by 
$$\left(0,\ e^{\theta(x,s)}\Lambda(dx, dt), \ \gamma(t)+\left(\int_0^t \int_{\mathbb{R}} (e^{\theta(x,s)}-1)x\Lambda(dx,ds)\right)'_t\right)_1$$
is equivalent to the distribution of $\xi$, i.e., 
$P_{\xi}\sim P_{\xi_{\theta}}$ and
the density transformation formula is of the form:
\begin{eqnarray*}
{\frac{d P_{\xi_{\theta}}}{dP_{\xi}}({\xi}(\cdot))}=\exp\left\{-\int_0^1\int_{{\mathbb{R}}\setminus \{0\}} \left(e^{\theta(x,s)}-1-\theta(x,s) \right) \Lambda(dx,dt)+\right. \\ \left. +\int_0^1 \int_{{\mathbb{R}}\setminus \{0\}} \theta(x,s) \bar{\cal{N}}_{\Lambda}(dx,dt)(\cdot) \right\},\ \ P_{\xi}\textit{ -a.e.}
\end{eqnarray*}   
\end{fact}
Comments: 
\begin{enumerate}
\item Condition (\ref{cond}) guarantees existence of the integrals and the properties of L\'evy measure for
$e^{\theta(x,s)}\Lambda(dx, dt)$.
\item Note that  if there exists $\theta^*$ such that
\begin{eqnarray}
\int_0^t \int_{{\mathbb{R}}\setminus \{0\}} (e^{\theta^*(x,s)}-1)x \Lambda(dx,ds)=-\int_0^t \gamma(s)\,ds,\quad \forall t\in[0,1],
\end{eqnarray}
then the transformed process is a martingale.
\end{enumerate}

\section{Shifted small ball probabilities for symmetric $\alpha$-stable processes}
\subsection{"Small" shifts.}\label{est}

Let $X_{\alpha}$ be a symmetric $\alpha$-stable L\'evy process, $\alpha\in (1,2)$. The aim of this section is to estimate shifted small ball
probabilities for these processes, unlike the centered small ball
probabilities that were studied in \cite{Mog_74}
\begin{eqnarray}\label{Mog}
\Pr \{X_{\alpha}\in \mathrm{B}(0,r) \}=\exp\{-{K_{\alpha}}{r^{-\alpha}}(1+o(1))\},
\end{eqnarray}
where $0<K_{\alpha}<\infty$, it depends just on the process
$X_{\alpha}$. This constant is equal to the first eigenvalue of
the fractional Laplacian (cf. \cite{ZRK_07}), the explicit
expression for $K_{\alpha} $ is still not found.
\begin{theorem}\label{small lambda}
For all $f\in \cal{C}$ and  $\lambda>0$, $r>0$ such that $ \lambda r^{\alpha-1}\to 0$, $ r\to 0$ we have
\begin{eqnarray*}
\Pr\left\{ \|X_{\alpha}(\cdot)-\lambda\, f(\cdot) \|< r \right\}=\exp\left\{-K_{\alpha}\,
r^{-\alpha}(1+o(1))\right\}.
\end{eqnarray*}
\end{theorem}
Comments:
\begin{enumerate}

\item In this theorem, we consider relatively small $\lambda$, namely $\lambda=o( r^{-(\alpha-1)})$. The case  when $\lambda$ is finite  is
included in this part of the result.

\item Notice that the estimate is similar with the estimate
(\ref{Mog}) for centered small deviations. The leading  term of the
asymptotic estimate is not sensitive for $f$, the dependence on
$f$ is hidden in the rest term.

\end{enumerate}

\begin{proof}

{\bf{Upper bound:}}  Using the Anderson inequality which holds for
symmetric processes (cf. \cite{LRZ_95}, \cite{BK_86}), and taking into account
(\ref{Mog}), we obtain
\begin{eqnarray*}
\Pr \{X_{\alpha}\in \mathrm{B}(\lambda\!\cdot\! f,r)\}\leq \Pr
\{X_{\alpha}\in \mathrm{B}(0,r)\}\leq
\exp\{-{K_{\alpha}}{r^{-\alpha}}(1+o(1))\}.
\end{eqnarray*}

{\bf{Lower bound:}}  We modify an approach from \cite{Shm_06}. 

Take $f\in \bf{H}$. By using
self-similarity, we can write
 \begin{eqnarray*}
 \Pr\{\| X_{\alpha}-\lambda\, f\|< r)\}= \Pr\{\| \xi_1\|< r\rho^{1/\alpha}\},
 \end{eqnarray*}
where $\xi_1$ is a L\'evy process with the centered
triplet $(0, \rho\, | x |^{-1-\alpha} dx\, dt,\
-\lambda\rho^{1/\alpha} f'(t))_1$ and $\rho$ is an arbitrary positive
real number that we are free to choose.

Using Fact~1, we have that an additive process $\xi_2$ with the generating triplet
\begin{eqnarray*}
 \left(0,\ \rho e^{\theta(x,t)}\frac{dx}{|x |^{1+\alpha} } dt, 0\right)_1,
\end{eqnarray*}
where $\theta(x,t)=\log \left(1+{\lambda \rho^{-\frac{\alpha-1}{\alpha} }}\cdot \frac{2-\alpha}{2}\,f'(t)\,x\, 1_{\{|x |<1\}} \right)$ has distribution $P_{\xi_2}$ equivalent to $P_{\xi_1}$ and it is a martingale. 
Note that it is defined correctly if
\begin{eqnarray*}
{\lambda \rho^{-\frac{\alpha-1}{\alpha} }}\cdot \frac{2-\alpha}{2}| f'(t)|<1 \textit{\ \ for almost all\ \ } t\in[0,1].
\end{eqnarray*}
If we assume
\begin{eqnarray}\label{*}
\rho: \   \lambda \rho^{-(\alpha-1)/\alpha }\to 0,
\end{eqnarray}
then this condition holds for large enough $\rho$.
This will be the first restriction we impose on $\rho$. 
We continue
\begin{eqnarray*}
 && \Pr\{X_{\alpha}-\lambda\, f\in B(0,r)\}=\int_{B(0,r\rho^{1/\alpha})} \frac{dP_{\xi_1}}{dP_{\xi_2}} dP_{\xi_2}= \\
 && =\exp\left\{-\rho\int_0^1 \int_{|\ell |<1}\Psi\left({\lambda \rho^{-\frac{\alpha-1}{\alpha} }}\cdot \frac{2-\alpha}{2}\,f'(t)\,x \right)\frac{d x}{|x |^{1+\alpha}}d t\right\}\times \\
 &&\E \exp\left\{-\int_0^1\int_{|x |<1} \log\left(1+{\lambda \rho^{-\frac{\alpha-1}{\alpha} }}\cdot \frac{2-\alpha}{2}\,f'(t)\,x\right)\,\bar{\cal{N}}_{\xi_2}(dx, dt)\right\} 1_{\{\| \xi_2\|< r\rho^{1/\alpha}\}}=\\
 &&= D_{\rho}\times S_{\rho}, 
\end{eqnarray*}
where $\Psi(u)= (1+u)\log(1+u)-u= \frac{u^2}{2} (1+o(1))$ as $u\to 0$. 

{\it{Deterministic term simplification.}}
\begin{eqnarray*}
  D_{\rho}=\exp\left\{- \frac{2-\alpha}{4} \lambda^2 \rho^{(2-\alpha)/\alpha} \int_0^1 f'(t)^2dt\ (1+o(1))\right\}.
\end{eqnarray*}

{\it{Stochastic term simplification.}}\ \  By the Jensen inequality, we get rid of the stochastic term of the density transformation formula
\begin{eqnarray*}
&& S_{\rho}\geq \exp\left\{-\E_{\Pr'} \int_0^1\int_{|x |<1} \log\left(1+{\lambda \rho^{-\frac{\alpha-1}{\alpha} }}\cdot \frac{2-\alpha}{2}\,f'(t)\,x\right)\,\bar{\cal{N}}_{\xi_2}(dx, dt)\right\}  \Pr \{{\xi_2} \in B(0, r \rho^{1/\alpha})\}=\\
&& = \Pr{\{\| \xi_2\|< r\rho^{1/\alpha}\}},
 \end{eqnarray*}
where $\Pr'$ : $ \frac{d\Pr'}{d\Pr}=1_{\{\|\xi_2\|<r\rho^{1/\alpha}\}}(\Pr{\{\|\xi_2\|<r\rho^{1/\alpha}\}})^{-1}$.
It is left just to treat the small ball probability $\Pr{\{\| \xi_2\|< r\rho^{1/\alpha}\}}$ for the time-inhomogeneous martingale $\xi_2$.

{\it{Homogenization.}}\quad
It is clear that ${\xi_2}(\cdot)\buildrel{d}\over{=} \xi(\rho\cdot)$, where $\xi$ is  a L\'evy process with the centered triplet
\begin{eqnarray*}
 \left(0,\ \left(1+{\lambda \rho^{-\frac{\alpha-1}{\alpha} }}\cdot \frac{2-\alpha}{2}\,f'(t)\,x\, 1_{\{|x |<1\}} \right)\frac{dx}{|x |^{1+\alpha} } dt, 0\right)_1.
\end{eqnarray*}
We can represent the process as a sum of  independent  processes
$\xi(\cdot)\buildrel{d}\over{=}\zeta_1(\cdot)+\zeta_2(\cdot),$  where $\zeta_1$ is a  L\'evy process generated
by the centered triplet
\begin{eqnarray*}
\left(0, \left(1- |x |  1_{\{|x|<1\}} \right)\frac{dx}{|x |^{1+\alpha}}, 0\right)_1,
\end{eqnarray*}
$\zeta_2$ is an additive process generated by
\begin{eqnarray*}
 \left(0,\ \left(|x | +{\lambda \rho^{-\frac{\alpha-1}{\alpha} }}\cdot \frac{2-\alpha}{2}\,f'(t)\,x\right) 1_{\{|x |<1\}} \frac{dx}{|x |^{1+\alpha} } dt, 0\right)_1.
\end{eqnarray*}
Taking into account this decomposition, for any $\delta\in (0, 1)$
we can write
\begin{eqnarray*}
 &\Pr \{\|{\xi(\rho\cdot)}\|< r \rho^{1/\alpha}\}\geq & \Pr \{\|\zeta_1(\rho\cdot)\|<(1-\delta) r \rho^{1/\alpha}\}
 \Pr \{\|\zeta_2(\rho\cdot)\|< \delta\cdot r \rho^{1/\alpha}\}.
\end{eqnarray*}
Let us treat each of the probabilities separately.

Using results of section~{8.2.4} of \cite{BGT_89}, one can prove that the process $\zeta_1$ belongs to the domain of normal attraction of
$X_{\alpha}$, i.e., 
\begin{eqnarray*}
\frac{\zeta_1(\rho \cdot)}{\rho^{1/\alpha}}\buildrel{d}\over{\Rightarrow} X_{\alpha}(\cdot) \textit{\ \ as\ \ } \rho\to \infty,
\end{eqnarray*}
where "$\buildrel{d}\over{\Rightarrow}$" means convergence in distribution. 

By \cite{Mog_74} and a slight generalization of his result
by \cite{Rus_07}, we know the following:

\begin{prop}\label{prop1}
For any L\'evy process $X$ which is a martingale and belongs the normal domain of attraction of a strictly
$\alpha$-stable L\'evy process $X_{\alpha}$, we have
\begin{eqnarray*}
\Pr\left\{ \frac{\|X(\rho \cdot)\|}{\rho^{1/\alpha}} < r
\right\} =
\exp\left\{-\frac{K_{\alpha}}{r^{\alpha}}(1+o(1))\right\},
\end{eqnarray*}
that holds  as $r \to 0$ and  $r \rho^{1/\alpha}\to \infty$. The constant
$K_{\alpha}$ is as in $(\ref{Mog})$.

In particular, if $X$ is from the normal domain of attraction to the
Wiener process, then $K_{2}={\pi^2}/{8}$.
\end{prop}

Thus, assuming
\begin{equation}\label{**}
\rho:\quad r \rho^{1/\alpha}\to \infty,
\end{equation}
we have
\begin{eqnarray*}
\Pr\left\{ \left\|\frac{\zeta_1 (\rho\cdot)}{\rho^{1/\alpha}}\right\| < (1-\delta)r \right\} = \exp\left\{-\frac{K_{\alpha}}{r^{\alpha}(1-\delta)^{\alpha} } (1+o(1))\right\}.
\end{eqnarray*}

In its turn, $\zeta_2$ could be decomposed into the sum of  processes with only  positive and only negative jumps, $\zeta_2
\buildrel{d}\over{=}\zeta^+ +\zeta^-$, where $\zeta^{\pm}$ are
generated by
\begin{eqnarray*}
\left(0,  \left(1 \pm \lambda \rho^{-\frac{\alpha-1}{\alpha} }\cdot \frac{2-\alpha}{2}\,f'(t)\right) 1_{\{\pm x \in (0,1)\}} \frac{dx}{(\pm x)^{\alpha}}  dt, 0 \right)_1.
\end{eqnarray*}
This decomposition give us
\begin{eqnarray*}
\Pr \{\|\zeta_2 (\rho\cdot)\|< \delta \cdot r \rho^{1/\alpha}\}\geq \Pr \{\|\zeta^+(\rho\cdot)\|< ({\delta}/{2}) r \rho^{1/\alpha}\}\Pr \{\|\zeta^-(\rho\cdot)\|< ({\delta}/{2}) r \rho^{1/\alpha}\}.
\end{eqnarray*}
Now, we need the following lemma, that allows  to switch to homogeneous processes.
\begin{lemma}\label{lemm}
Let $\Lambda$ be a L\'evy measure such that $\int_{\vert x \vert>1} x \Lambda(dx)<\infty$ and $\mu(\cdot)\in AC[0,1]$.

If $\eta$ is an additive process specified by  the generating triplet $(0, \mu'(t) dt\, \Lambda(d x), 0)_1$ and
$\zeta$ is a L\'evy process specified by the generating triplet $(0, (\int_0^1\mu(t) dt) \,\Lambda(d x), 0)_1$, then
\begin{eqnarray*}
\|\eta\|\buildrel{d}\over{=}\|\zeta\|.
\end{eqnarray*}
\end{lemma}
\begin{proof}
By L\'evy-Khintchine formula  we have
\begin{eqnarray*}
\|\eta\|\buildrel{d}\over{=} \sup_{t\in[0,1]}\left|X\left(\int_0^t\mu(s)ds\right)\right|=\sup_{t\in [0,\int_0^1\mu(s)ds]}\left|X\left(\frac{\int_0^t\mu(s)ds}{\int_0^1\mu(s)ds}\right)\right|,
\end{eqnarray*}
where $X$ is a L\'evy process generated by $(0,
\Lambda(d x), 0)_1$.

Put $\varphi(t)={\int_0^t\mu(s)ds}/{\int_0^1\mu(s)ds}$. Notice that
$0\leq \varphi(t)\leq 1$ for all $t\in [0,1]$. Taking this into account
we continue
\begin{eqnarray*}
\|\eta\|\buildrel{d}\over{=}\sup_{t\in[0,\int_0^1\mu(s)ds]}|X(\varphi(t))|\buildrel{d}\over{=} \sup_{t\in[0,1]}|\zeta(\varphi(t))|=\sup_{s\in[0,1]}|\zeta(s)|.
\end{eqnarray*}
\end{proof}
Using this lemma, we continue
\begin{eqnarray*}
\Pr \{\|\zeta_2 (\rho\cdot)\|< \delta \cdot r \rho^{1/\alpha}\}\geq \Pr \{\|\eta^+(\rho\cdot)\|< ({\delta}/{2}) r \rho^{1/\alpha}\}\Pr \{\|\eta^-(\rho\cdot)\|< ({\delta}/{2}) r \rho^{1/\alpha}\},
\end{eqnarray*}
 where $\eta^{\pm}$
are centered subordinators generated by
\begin{eqnarray*}
\left(0,  \left(1 \pm \lambda \rho^{-\frac{\alpha-1}{\alpha} }\cdot \frac{2-\alpha}{2}\,f(1)\right) 1_{\{\pm x \in (0,1)\}} \frac{dx}{(\pm x)^{\alpha}}  dt, 0 \right)_1.
\end{eqnarray*}
The processes $\eta^{\pm}$ have just bounded jumps,
therefore they both belong to the normal domain of attraction of the
Wiener process. By Proposition \ref{prop1}, under (\ref{**}) we obtain
\begin{eqnarray*}
\Pr \{\|\eta^{\pm}(\rho\cdot)\|< (\delta/2) r \rho^{1/\alpha}\}\geq \exp\left\{-\frac{\pi^2}{2\delta^2(3-\alpha)}
\, r^{-2} \rho^{-\frac{2-\alpha}{\alpha}}\left(1 \pm \lambda \rho^{-\frac{\alpha-1}{\alpha} }\cdot \frac{2-\alpha}{2}\,f(1)\right)(1+o(1))\right\} .
\end{eqnarray*}
Then
\begin{eqnarray*}
\Pr \{\|\zeta_2 (\rho\cdot)\|< \delta \cdot r \rho^{1/\alpha}\} \geq \exp\left\{-\frac{\pi^2}{\delta^2(3-\alpha)}
\cdot \frac{1}{r^2 \rho^{(2-\alpha)/\alpha}}(1+o(1))\right\}.
\end{eqnarray*}
Thus, we obtain that  for any $\delta\in (0,1)$
\begin{eqnarray*}
\Pr \{\|\xi_2\|< r \rho^{1/\alpha}\} \geq \exp\left\{
-\frac{K_{\alpha}}{r^{\alpha}(1-\delta)^{\alpha} } (1+o(1))-\frac{\pi^2}{\delta^2(3-\alpha)}
\cdot \frac{1}{r^2 \rho^{(2-\alpha)/\alpha}}(1+o(1))\right\}.
\end{eqnarray*}
Under (\ref{**}) we have
\begin{eqnarray*}
\Pr \{\|\xi_2\|< r \rho^{1/\alpha}\} \geq \exp\left\{
-\frac{K_{\alpha}}{r^{\alpha}} (1+o(1))\right\}.
\end{eqnarray*}
Collecting all the preliminary results, under (\ref{*}) and (\ref{**}) we get
\begin{eqnarray*}
  \Pr\left\{ \|X_{\alpha}(\cdot)-\lambda\, f(\cdot) \|< r \right\}\geq \exp\left\{- \frac{2-\alpha}{4} \lambda^2 \rho^{(2-\alpha)/\alpha} \int_0^1 f'(t)^2dt\ (1+o(1))-\frac{K_{\alpha}}{r^{\alpha}} (1+o(1))\right\}.
\end{eqnarray*}
Using the condition $\lambda r^{\alpha-1}\to 0$ of the theorem  (for the first time in the proof), we can find 
 $\rho$ obeying (\ref{*}) and (\ref{**}), such that  $\lambda^2 \rho^{(2-\alpha)/\alpha}=o(r^{-\alpha})$. For example, take $\rho^*= r^{-\alpha}(\lambda r^{\alpha-1})^{-1}$. 
Thus,  for all $f\in \bf{H}$ we have
\begin{eqnarray*}
  \Pr\left\{ \|X_{\alpha}(\cdot)-\lambda\, f(\cdot) \|< r \right\}\geq \exp\left\{-\frac{K_{\alpha}}{r^{\alpha}} (1+o(1))\right\}.
\end{eqnarray*}
The set $\bf{H}$ is dense in $\cal{C}$, so the result could be generalized for arbitrary $f\in \cal{C}$.
\end{proof}
Remark:

We can also exploit the same proof 
under $ \lambda r^{\alpha-1}\to \infty$ or  $ \lambda r^{\alpha-1}=c, \ c>0$ conditions.

For example, under $ \lambda r^{\alpha-1}\to \infty$ for $f'\in L_{\infty}$
taking $\rho:=\left(\lambda\|f'\| (2-\alpha)
2(1-\epsilon)\right)^{\alpha/(\alpha-1)}$,  $\epsilon\in(0,1)$,  we get
\begin{eqnarray*}
\Pr\left\{ \|X_{\alpha}(\cdot)-\lambda\, f(\cdot) \|< r \right\}\geq
 \exp\left\{- C_1 \, \lambda^{\alpha/(\alpha-1)}(1+o(1))\right\},
\end{eqnarray*}
where
\begin{eqnarray*}
C_1=C_1(f, \alpha) = \|f'\|^{\alpha/(\alpha-1)}\left(\frac{2-\alpha}{2}\right)^{\alpha/(\alpha-1)}
\sum_{k=1}^{\infty}\frac{\int_0^1( f'(t)/\|f'\|)^{2k}dt}{k(2k-1)(2k-\alpha)}.
\end{eqnarray*}
We see that the order differs from the order of the upper estimate, and moreover, we can prove that it is not optimal.  Following  \cite{AD_08} we can obtain: 
there are constants $0<C_1 \leq C_2<\infty $ s.t.
\begin{eqnarray*}
\exp\left\{- C_2 \cdot \frac{\lambda}{r}\log \lambda r ^{\alpha-1}\right\}\leq\Pr\left\{ \|X_{\alpha}(\cdot)-\lambda\, Id(\cdot) \|< r \right\}\leq
 \exp\left\{- C_1 \cdot \frac{\lambda}{r}\log \lambda r ^{\alpha-1}\right\},
\end{eqnarray*}
where $Id$ is the identity function on $[0,1]$. 

Under $ \lambda r^{\alpha-1}=c, \ c>0$ condition, we are  faced with a known open problem for processes  from the  domain of attraction of  $X_{\alpha}$
\begin{eqnarray*}
\Pr\left\{ {\|X(\rho \cdot)\|}< c
\right\} =
\exp\left\{- A_{\alpha}(c) \rho(1+o(1))\right\},
\end{eqnarray*}
$A_{\alpha}(c)$ is not known here.

Nevertheless, in the next section we obtain a result in this case, slightly modifying the proof.

\subsection{"Middle" shifts.}
\begin{theorem}\label{LIL case}
For any $c>0$, $f\in AC[0,1]\ :\  f'\in L_{\infty},$ $f(0)=0$ such that
\begin{eqnarray}\label{ref}
\|f'(\cdot)\|<\frac{2}{2-\alpha}\cdot \frac{1}{c},
\end{eqnarray}
we have
\begin{eqnarray*}
\exp\left\{- K_{\alpha}
\,\frac{1}{r^{\alpha}}(1+o(1))\right\}\geq\Pr\left\{ \|X_{\alpha}(\cdot)- c\cdot r^{-(\alpha-1)}\, f(\cdot)
\|< r \right\}\geq \exp\left\{- C(\alpha)
\,\frac{1}{r^{\alpha}}\right\},
\end{eqnarray*}
as $ r\to 0$, where
\begin{eqnarray*}
C(\alpha)= 2\left(\frac{1}{\alpha}+ \sum_{k=1}^{\infty} \frac{1}{2k(2k-1)(2k-\alpha)}+24 \cdot 6^{\alpha}\left(\frac{1}{2-\alpha}+\frac{2^{\alpha-1}-1}{6(3-\alpha)}\right)\right),
\end{eqnarray*}
and $K_{\alpha}$ is as in $(\ref{Mog})$.
\end{theorem} 
\begin{proof}
{\bf{Upper bound:}}  The Anderson inequality.

{\bf{Lower bound:}}  In this proof, we are close to \cite{AD_08} method.  
We start with a truncation of large jumps
\begin{eqnarray*}
\Pr\left\{ \|X_{\alpha}(\cdot)- c\cdot r^{-(\alpha-1)}\, f(\cdot)
\|< r \right\}= \Pr\{\|X_{\alpha}(\cdot)- c\cdot r^{-(\alpha-1)}\, f(\cdot)
\|< r\, \vert\ A\}\Pr\{A\},
\end{eqnarray*}
where $A$ is the event that the process $X_{\alpha}$ has no jumps bigger than $r$, i.e., $A=\{\omega\in\Omega :  \forall  t \in[0,1]\  \Delta X_{\alpha}(t,\omega)\leq r\}$.
It is well-known that
\begin{eqnarray*}
\Pr\{A\}=\exp\left\{- \int_{\vert x \vert>r} \frac{dx}{\vert x\vert^{\alpha+1}}\right\}= \exp\left\{-\frac{2}{\alpha}\cdot r^{-\alpha}\right\}.
\end{eqnarray*}
Denote by $\xi_1$ an additive process with the generating triplet $$(0,\  1_{\{\vert x \vert< r\}} \vert x \vert^{-(1+\alpha)} dx\, dt ,\  - c\, r^{-(\alpha-1)}f'(t))_1.$$
Hence, we continue 
\begin{eqnarray*}
\Pr\left\{ \|X_{\alpha}(\cdot)- c\cdot r^{-(\alpha-1)}\, f(\cdot)
\|< r \right\}=\exp\left\{-\frac{2}{\alpha}\cdot r^{-\alpha}\right\} \Pr\{\|\xi_1\|< r\}.
\end{eqnarray*}
Using Fact~1, we obtain that $\xi_2$  with the generating triplet 
$$\left(0,\  \left(1+c\cdot \frac{2-\alpha}{2} f'(t)\frac{x}{r} \right)1_{\{\vert x \vert< r\}} \frac{dx}{\vert x \vert^{1+\alpha}} dt ,\  0\right)_1$$
has  distribution equivalent to $P_{\xi_1}$. Take $\theta^*(s,t)=\log(1+c\cdot \frac{2-\alpha}{2} f'(t)\frac{x}{r})$ in Fact ~1.
Note that $\xi_2$ is correctly defined under condition (\ref{ref}), and $\xi_2$ is a martingale.

Thus, we continue
\begin{eqnarray*}
&&\Pr\{\|\xi_1\|<r\}= \int_{B(0,r)} \frac{dP_{\xi_1}}{dP_{\xi_2}}(\eta) dP_{\xi_2}(\eta)=\\
&&=\exp\left\{-\int_0^1\int_{\vert x \vert<r} \Psi\left(c\cdot \frac{2-\alpha}{2} f'(t)\frac{x}{r}\right)\frac{dx}{\vert x \vert^{1+\alpha}} dt \right\} \times \\
&& \times \, \E \exp\left\{-\int_0^1\int_{\vert x \vert<r} \log\left(1+c\cdot \frac{2-\alpha}{2} f'(t)\frac{x}{r}\right)\bar{N}_{\Lambda_2}(dx, dt) \right\} 1_{\{\|\xi_2\|<r\}}=\\
&&=D \times S,
\end{eqnarray*}
where $\Psi(u)=(1+u)\log(1+u)-u=\sum_{k=2}^{\infty}(-1)^k\frac{u^k}{k(k-1)}$, as $\vert u\vert<1$.

{\it{Deterministic term simplification.}}
\begin{eqnarray*}
&& D=\exp\left\{-\frac{2}{r^{\alpha}} \sum_{k=1}^{\infty}\int_0^1\left(c\cdot \frac{(2-\alpha)  f'(s) }{2}\right)^{2k}ds\cdot  \frac{1}{2k(2k-1)(2k-\alpha)} \right\}\geq\\
&& \geq \exp\left\{-\frac{2}{r^{\alpha}} \sum_{k=1}^{\infty} \frac{1}{2k(2k-1)(2k-\alpha)} \right\}.
\end{eqnarray*}

{\it{Stochastic term simplification.}}
Take  a probability measure $ \Pr' : \  \frac{d\Pr'}{d\Pr}=\frac{1_{\{\|\xi_2\|<r\}}}{\Pr{\{\|\xi_2\|<r\}}}$, then by Jensen's inequality obtain
\begin{eqnarray*}
S\geq \exp\left\{-\E_{\Pr'} \int_0^1\int_{\vert x \vert<r} \log\left(1+c\cdot \frac{2-\alpha}{2} f'(t)\frac{x}{r}\right)\bar{N}_{\Lambda_2}(dx, dt) \right\} \Pr{\{\|\xi_2\|<r\}}= \Pr{\{\|\xi_2\|<r\}}.
\end{eqnarray*}
To estimate the last probability, we use a proposition proved in \cite{AD_08} (see Lemma 4.1 there)
\begin{prop}\label{martingale}
 Let  X be a $(\nu, 0)$-L\'evy martingale with $\nu$ supported on $[-\varepsilon, \varepsilon]$,  then
\begin{eqnarray*}
\Pr{\{\| X\|<3\varepsilon\}}\geq \exp\left\{-\left(12\, \frac{1}{\varepsilon^2}\int_{|x|<\varepsilon} x^2 \nu(dx)+2\right)\right\}.
\end{eqnarray*}
\end{prop}
We can't use the proposition directly, because $\xi_2$ being a martingale nevertheless is time-inhomogeneous.  

{\it{Homogenization.}} We decompose the process $\xi_2$ into a sum of independent processes,  one of which $\zeta_1$ is a L\'evy process with the L\'evy measure
\begin{eqnarray*}
\left(1- \frac{|x|}{r}\right)1_{\{\vert x \vert< r\}} \frac{d x}{|x |^{1+\alpha}},
\end{eqnarray*}
 and the second $\zeta_2$ is an additive process with the L\'evy measure 
 \begin{eqnarray*}
 \left(|x| + c\cdot \frac{2-\alpha}{2}f'(t) x \right)\frac{1}{r}1_{\{\vert x \vert< r\}} \frac{d x}{|x |^{1+\alpha}}dt,
 \end{eqnarray*}
 we can choose shifts in such a way that both of the processes are again martingales
\begin{eqnarray*}
\xi_2(\cdot)\buildrel{d}\over{=}\zeta_1(\cdot)+\zeta_2(\cdot).
\end{eqnarray*}
Taking into account this decomposition, for any $\delta\in (0, 1)$
we can write
\begin{eqnarray*}
\Pr\{\|\xi_2\|< r \} \geq  \Pr \{\|\zeta_1\|<(1-\delta) r\} \Pr \{\|\zeta_2\|< \delta r\}.
\end{eqnarray*}

 Let us treat each of the probabilities separately.

Using Proposition \ref{martingale}, we obtain for $\zeta_1$
 \begin{eqnarray*}
\Pr \{\|\zeta_1\|<(1-\delta) r\}\geq\exp\left\{ -\frac{1}{r^{\alpha}}\cdot 24 \left(\frac{3}{1-\delta}\right)^{\alpha} \left(\frac{1}{2-\alpha}-\frac{1-\delta}{3(3-\alpha)}\right)\right\}.
\end{eqnarray*}

The sample paths of $\zeta_2$ are of bounded variation. Thus, the following decomposition into the sum of  processes with only positive and only negative jumps is possible
 $\zeta_{2}\buildrel{d}\over{=}\zeta^+ +\zeta^-$, where $\zeta^{\pm}$ are
generated by
\begin{eqnarray*}
\left(0,\ \left(1\pm c\cdot \frac{2-\alpha}{2} f'(t) \right)\frac{1}{r} \, 1_{\{\pm x \in (0, r)\}} \frac{dx}{(\pm x)^{\alpha}} dt,\  0\right)_1,
\end{eqnarray*}
correspondingly. This decomposition yields
\begin{eqnarray*}
\Pr \{\|\zeta_2\|< \delta r \}\geq \Pr \{\|\zeta^+\|< ({\delta}/{2}) r \}\Pr \{\|\zeta^-\|< ({\delta}/{2}) r \}.
\end{eqnarray*}
Using Lemma \ref{lemm}, we continue
\begin{eqnarray*}
\Pr \{\|\zeta^{\pm}\|< ({\delta}/{2}) r \}=\Pr \{\|\eta^{\pm}\|< ({\delta}/{2}) r \},
\end{eqnarray*}
 where $\eta^{\pm}$ are centered  positive (negative) subordinators generated by
\begin{eqnarray*}
\left(0,\ \left(1\pm c\cdot \frac{2-\alpha}{2} f(1) \right)\frac{1}{r} \, 1_{\{\pm x \in (0, r)\}} \frac{dx}{(\pm x)^{\alpha}},\ 0\right)_1.
\end{eqnarray*}
Applying Proposition \ref{martingale}, obtain 
\begin{eqnarray*}
\Pr \{\|\zeta_2\|< \delta r \}\geq \exp\left\{-\frac{1}{r^{\alpha}} \cdot  \frac{24}{3-\alpha} \cdot \left(\frac{6}{\delta}\right)^{\alpha-1} \right\}.
\end{eqnarray*}
For simplicity, take $\delta=1/2$ and obtain the statement of the theorem.
\end{proof}

 \section{Law of the Iterated Logarithm for stable L\'evy processes}\label{LIL}
 \subsection{General information.}

There are several recent works that deal with non-standard Law of
the Iterated Logarithm  (LIL) statements for L\'evy processes and
random walks, in particular, in the case when the variance of
random variables is infinite, see \cite{Ein_07}, \cite{BDM_08},
\cite{Sav_08}, \cite{CKL_00}.

In this section, we collect facts related to the LIL for the stable
L\'evy processes.   
Traditionally LIL statements could be of 
Limsup (Strassen) or Liminf (Chung) types.  Proofs of the first type of results are  based on large deviation inequalities, whereas the second type of results usually needs small deviation estimates. 

{\it{Limsup LIL:}}

One of the interpretations of the LIL is
the rate of convergence in the CLT theorem. Analogue of the
functional  CLT theorem (invariance principle) for stable
processes is:
\begin{eqnarray*}
\frac{X(T\cdot)}{T^{1/\alpha}L(T)}\buildrel{d}\over{\Rightarrow}
X_{\alpha}(\cdot),
\end{eqnarray*}
where $X$ is a process form the domain of attraction of $X_{\alpha}$ and
$L(\cdot)$ is a proper slowly varying function. If $X$ is $X_{\alpha}$ itself, this relation is nothing more than the  self-similarity property
\begin{eqnarray*}
\frac{X_{\alpha}(T\cdot)}{T^{1/\alpha}}\buildrel{d}\over{=}
X_{\alpha}(\cdot).
\end{eqnarray*}
The Marzinkevich-Zygmund LLN says
\begin{eqnarray*}
\frac{X_{\alpha}(T)}{T^{1/p}}{\rightarrow}\ 0 \textit{\ \ \ a.s.,\
\ \ if \ \ } p\in (1,\alpha).
\end{eqnarray*}
As for the LIL, the situation is predetermined by the following
dichotomy statement (cf. Thoerem VIII.5 in  \cite{Ber_96}):
\begin{fact}\label{fact 2}
\begin{eqnarray*}
\limsup_{T\to \infty}\frac{|X_{\alpha}(T)|}{T^{1/\alpha}h(T)}=0 \textit{ \ a.s. \qquad or }\ =\infty \textit{\ \ a.s.}
\end{eqnarray*}
according as
 \begin{eqnarray*}
\int^{\infty} \frac{d\,t}{t h(t)^{\alpha}}<\infty \textit{ \qquad or }\ =\infty.
\end{eqnarray*}
\end{fact}

This fact says that the stable L\'evy processes doesn't exhibit
LIL behavior: there is no such a function $\varphi(\cdot)$ that
$0<\limsup_{t\to\infty}{|X_{\alpha}(T)|}/{\varphi(T)}<\infty$.


The statemet gives the following information on the
sample paths growth at infinity: according to the integral test
\begin{eqnarray*}
\Pr\{\omega : \exists t_0(\omega)\textit{\ s.t.\  for}  \  t\in(t_0(\omega),
\infty)\ \ |X_{\alpha}(t,\omega)|< t^{1/\alpha}h(t)\}=1 \textit{\ \ or\ \ } =0.
\end{eqnarray*}
That also means that the set $ \{t:\ |X_{\alpha}(t)|>
t^{1/\alpha}h(t)\}$ is a.s. bounded or unbounded according to the
integral  test.

For example, we can say that almost all sample paths of the process  $X_{\alpha}(t), t\in (0,\infty)$ intersect the level $\varphi(t)=t^{1/\alpha}(\log t)^{1/\alpha}$
infinitely many times, whereas the level $\psi(t)=t^{1/\alpha}(\log t)^{\epsilon +1/\alpha}$, $\epsilon>0$ is overpassed just finitely many times.

In what follows, we need a limsup statement for the sup-process $M(\cdot)$ that is  an increasing sample paths process defined by
 $$
 M(T)={\sup_{s\in [0,1]}| X_{\alpha}(T s)|}.
 $$ 
\begin{corol}
For any $\varphi: \  \mathbf{R}^+ \to \mathbf{R}^+$ s.t. $\int^{\infty}{dx}/{\varphi(x)}=\infty$ the following holds  
\begin{equation}\label{up_norm}
\limsup_{T\to\infty} \frac{ M(T) }{T^{1/\alpha}(\log T \cdot \varphi(\log\log T))^{1/\alpha}}=\infty \textit{\ \ a.s.}
\end{equation} 
\end{corol}

{\it{Liminf LIL:}}
  
Despite the fact that the standard LIL doesn't exist, the Chung-type LIL  for the stable L\'evy processes holds
\begin{eqnarray}\label{taylor}
\liminf_{T\to \infty}\frac{M(T) }{(T/\log\log T)^{1/\alpha} }=K_{\alpha}^{1/\alpha} \textit{\ \ a.s.,\ \ }
\end{eqnarray}
where $K_{\alpha}$ is as in (\ref{Mog}).
The law was discovered in \cite{Tay_67}.  This statement is about the rate of  moving of the sup-process away from zero. More precisely,  almost all sample paths of the sup-process finitely often intersect the level 
 $(1-c)K_{\alpha}^{1/\alpha}({T}/{\log\log T})^{1/{\alpha}}$
 and infinitely often $(1+c)K_{\alpha}^{1/\alpha}({T}/{\log\log T})^{1/{\alpha}}$, for any $0<c<1$, i.e.,
\begin{eqnarray*}
\Pr\{\omega:\  \{T:  M(T, \omega)<(1-c)K_{\alpha}^{1/\alpha}({T}/{\log\log T})^{1/{\alpha}}\}\ \ bdd \}=1,\\
 \Pr\{\omega:\  \{T:  M(T, \omega)<(1+c)K_{\alpha}^{1/\alpha}({T}/{\log\log T})^{1/{\alpha}}\}\ \ unbdd \}=1.
 \end{eqnarray*}
Combining (\ref{up_norm}) Êand (\ref{taylor}), we can say that for any $c\in (0,1)$, any $\varphi$ s.t. $\int^{\infty}{dx}/{\varphi(x)}=\infty$ the following holds:
for T large enough 
 \begin{eqnarray*}
M(T)\in \left((1-c)K_{\alpha}^{1/\alpha}(T/{\log\log T})^{1/{\alpha}},\quad T^{1/\alpha}(\log T)^{1/\alpha}(\varphi(\log\log T))^{1/\alpha}\right) \ \ \  a.s.
\end{eqnarray*}

In this article, we study a generalization of these results to a functional LIL. 
What we get is  analogous  to the result of Baldi and Royonette in Gaussian case \cite{BR_92}.

{\it{Baldi-Royonette  result for the Wiener process:}}
By  $W$ denote the Wiener process. Consider a family of scaling of $W$
\begin{eqnarray*}
\xi_T^{\gamma}(\cdot)= \frac{W(T\cdot)}{\sqrt{2T \log\log T}}\cdot \gamma(T),
\end{eqnarray*}
where $\gamma:  \  \mathbf{R}^+ \to \mathbf{R}^+$ s.t. $\gamma(0)=0$.  

{\bf{Definition:}}\quad Let  $(E,\tau)$  be a topological space. 
An element $x\in E$ is called an a.s. limit point of a family $\{\xi_T\}_{T>0}$ of random elements on $E$,  if  there exists $\{T_k\}_{k=1}^{\infty}$, $T_k\to \infty$ such that  $\xi_{T_k}(\omega)\buildrel{\tau}\over{\to} x$ as $k \to \infty$,  for almost all $\omega\in\Omega$.\\
The set of all a.s. limit points  of $\{\xi_T\}_{T>0}$, say ${\cal{K}}$, is called
the a.s. limit (cluster) set of   $\{\xi_T\}_{T>0}$. We  write $\{\xi_T\}_{T>0}\to\to {\cal{K}}$.

If we deal with $C[0,1]$ endowed with the uniform topology (it is known to be separable), then $\{\xi_T\}_{T>0}\to\to {\cal{K}}$ iff
\begin{enumerate}
\item $\lim_{T\to \infty} \inf_{f\in {\cal{K}}} \|\xi_T-f\|=0$ a.s.,  and
\item 
for all $f\in {\cal{K}}$\ \
$\liminf_{T\to \infty}\|\xi_T- f\|=0$ a.s.
\end{enumerate}

Depending on the rate of growth of $\gamma(\cdot)$ the following variants of a.s. cluster sets for the family 
$\{\xi_T^{\gamma}\}_{T>0}$ exist:
\begin{description}
\item[ (a)] If $\gamma(T)=o(1)$, then the cluster set  consists just from the zero function, which we denote by  $\mathbf{0}$
$$\{\xi_T^{\gamma}\}\to\to \{\mathbf{0}\}.$$
\item[ (b)] If $\gamma(T)\to c$, then the a.s. cluster set is a compact. Namely,  the Strassen LIL  holds
$$\{\xi_T^{\gamma}\}\to\to c^2 \, {\cal{S}},$$
where ${\cal{S}}= \{f\in{\bf{H}}, \int_0^1f'^2\leq 1\}.$
\item[ (c)] If $\gamma(T)\to \infty$ in such a way that  $\gamma(T)=o(\log \log T)$, then the following is true:  
$$\{\xi_T^{\gamma}\}\to\to \cal{C}.$$
\item[ (d)]
If $\gamma(T)\to \infty$ in such a way that there is $c_0>0$  such that $ c_0 \log \log T\leq \gamma(T)$ for large enough T, then the cluster set is empty. Namely, for any $f\in \cal{C}$ we have
\begin{eqnarray*}
 \liminf_{T\to\infty}\, \|\xi_T^{\gamma}(\cdot)-f(\cdot)\| \geq  \frac{c_0 \pi}{4}\ \ a.s.
\end{eqnarray*}
This scaling is too small to overpower natural fluctuations of the Wiener process, that is why the trajectories stop a.s. clustering around continuous functions. 
\end{description}

\subsection{Functional LIL for scaled stable L\'evy processes.}
In this section, we work in $D[0,1]$ endowed with the uniform topology, this is known non-separable topological space. The process $X_{\alpha}$ has no time-fixed jumps, therefore the uniform convergence is possible just to continuous functions. In this case, a.s. cluster sets ${\cal{K}}_h=\{f\in D[0,1]: \ \liminf_{T\to\infty}\|\frac{X_{\alpha}(T\cdot)}{T^{1/\alpha}h(T)} -f(\cdot)\|=0 \}$, if exist, are contained in $\cal{C}$. 

\begin{theorem}\label{Chung_deg}
Let $\  h: \ \mathbf{R}^+ \to \mathbf{R}^+$ such that $h(0)=0$ and there exists $c>0$ such that $h(T)\leq c (\log\log T)^{-1/\alpha}$ for $T$ large enough.
Then,  for any $f\in \cal{C}$ the following holds
\begin{eqnarray*}
\liminf_{T\to\infty}\, \left\|\frac{X_{\alpha}(T\cdot)}{T^{1/\alpha}h(T)} -f(\cdot)\right\|\geq \frac{K_{\alpha}^{1/\alpha}}{c} \textit{\ a.s.},
\end{eqnarray*}
where $K_{\alpha}$ is as in (\ref{Mog}).
\end{theorem}

This statement  corresponds to the case $(d)$ for the Wiener process. This is a degenerate situation from the point of view of the functional LIL, fluctuations of the process overpower the scaling, it corresponds to the empty limit set. 

If $h(\cdot)$ is such that $\lim_{T\to\infty}\frac{|X_{\alpha}(T)|}{T^{1/\alpha}h(T)}=0$ in Fact 2, then the a.s. cluster set is not bigger than $\{\mathbf{0}\}$.  It corresponds to the case $(a)$ for the Wiener process (scaling is too strong).

We eliminate these two well-understood cases. Hence,  our interest is focused  on the set of scaling functions  
$h: \ \mathbf{R}^+ \to \mathbf{R}^+$  obeying conditions:
$h(0)=0$ and for any $c>0$, any $\varphi:\ \int^{\infty}{dx}/{\varphi(x)}=\infty$  there exists $t_0>0$ such that  for all $T\in(t_0, \infty)$
 \begin{eqnarray*}
 c (\log\log T)^{-1/\alpha}<h(T)\leq  (\log T)^{1/\alpha}(\varphi(\log\log T))^{1/\alpha} .
 \end{eqnarray*}
\\
Small deviations estimates from Theorem \ref{small lambda} give us the following statement.
\begin{theorem}\label{LILL}
 For any $ f \in \cal{C}$, any 
$\delta\in (0,1)$ we have
\begin{eqnarray*}
\liminf_{T\to\infty}\, (\log\log T)^{\delta} \left\|\frac{X_{\alpha}(T\cdot)}{T^{1/\alpha}( \log\log T)^{\delta-1/\alpha}} -f(\cdot)\right\|= K_{\alpha}^{1/\alpha} \textit{\ a.s.}
\end{eqnarray*}
\end{theorem}
Comments: \begin{enumerate}
\item From this statement it follows that if $\delta\in (0,1),$ then 
$$\left\{\frac{X_{\alpha}(T\cdot)}{T^{1/\alpha}( \log\log T)^{\delta-1/\alpha}}\right\} \to\to\cal{C}.$$ 

\item Take  $\delta=1/\alpha$, to obtain the following effect:
 $$\left\{\frac{X_{\alpha}(T\cdot)}{T^{1/\alpha}}\right\} \to\to\cal{C},$$
despite  the fact that $\frac{X_{\alpha}(T\cdot)}{T^{1/\alpha}}\buildrel{d}\over{=} X_{\alpha}(\cdot)$ for any $T>0$. 
Under any fixed $T$ we get c\`adl\`ag looking trajectories, nevertheless the gradual scaling (moving $T$ to $\infty$)  of  the trajectories causes their clustering (almost all of them) around  continuous functions. The same effect took place  for $W$, see (c) case  under $\gamma(T)=\sqrt{\log\log T}$.

\item We already mention that the uniform topology is not separable on $D[0,1]$, therefore the a.s. cluster set could be bigger if we consider the Skorokhod topology, which is separable on $D[0,1]$. But anyway the cluster set will contain $\cal{C}$ because the uniform convergence implies convergence in the Skorokhod topology on $D[0,1]$.    
\end{enumerate}

Small deviations estimates from Theorem \ref{LIL case} give us the following statement
\begin{theorem}\label{LIL}
For any $f$ that belongs to
 \begin{eqnarray*}
{\cal{C}^*}=\left\{ f\in AC[0,1]:\ f(0)=0,\  \|f'\|<\frac{2}{2-\alpha}\cdot (C(\alpha))^{-(\alpha-1)/\alpha}\right\}
\end{eqnarray*}
we have
\begin{eqnarray*}
\liminf_{T\to\infty}\, (\log\log T) \left\|\frac{X_{\alpha}(T\cdot)}{(\,T/\log\log T)^{1/\alpha} \log\log T} -f(\cdot)\right\|= C'\ \ a.s.,
\end{eqnarray*}
where $C'\in [K_{\alpha}^{1/\alpha},(C(\alpha))^{1/\alpha}]$, $C(\alpha)$ is from Theorem \ref{LIL case} and $K_{\alpha}$ is as in $(\ref{Mog})$.
\end{theorem}
Comment:\\  From this statement it follows that the a.s. limit set of $\left\{\frac{X_{\alpha}(T\cdot)}{T^{1/\alpha}( \log\log T)^{1-1/\alpha}}\right\}_{T>0}$ contains ${\cal{C}^*}$.

\begin{proof}   We modify the proofs of Theorem VIII.6 in \cite{Ber_96} and Theorem 17.1 in \cite{Lif};  for the lower bound we also use ideas of \cite{Csa_80}.

{\bf{Lower bound in Theorems \ref{Chung_deg},  \ref{LILL} and \ref{LIL}:}}  
 Let $\delta\in [0,1]$, where $\delta=0$ corresponds to Theorem \ref{Chung_deg}, $\delta\in(0, 1)$  to Theorem \ref{LILL}  and $\delta=1$ to Theorem \ref{LIL}. Choose  $T_k= \exp\{k (\log k)^{-3}\}$.  We start with the inequalities
\begin{eqnarray*}
&& \liminf_{T\to\infty}\, (\log\log T)^{\delta} \left\|\frac{X_{\alpha}(T\cdot)}{T^{1/\alpha} (\log\log T)^{\delta-\frac{1}{\alpha}}} - f(\cdot)\right\| \geq\\
&& \liminf_{k\to\infty}\,\frac{ \inf_{[T_k, T_{k+1}]} \left\|X_{\alpha}(T\cdot) -f(\cdot)\,T^{1/\alpha} (\log\log T)^{\delta-\frac{1}{\alpha}}\right\|}{(T_{k+1}/\log\log T_{k+1})^{1/\alpha}}\geq\\
&& \liminf_{k\to\infty}\,\frac{ \left\|X_{\alpha}(T_k\cdot) -f(\cdot)T_k^{1/\alpha} (\log\log T_k)^{\delta-1/\alpha}\right\|}{(T_{k+1}/\log\log T_{k+1})^{1/\alpha}}.
\end{eqnarray*}
For the last inequality we used the following technical lemma
\begin{lemma}\label{lemma}
\begin{itemize}
\item[$(i)$]
For any increasing sequence $\{T_k\}_{k>0}$, any $f\in \cal{C}$  there exists $0<M<\infty$ such that the following is true
\begin{eqnarray*}
&& \inf_{[T_k, T_{k+1}]} \left\|X_{\alpha}(T\cdot) - f(\cdot)T^{1/\alpha} (\log\log T)^{\delta-1/\alpha}\right\|\geq
\left\|X_{\alpha}(T_k\cdot) - f(\cdot)T_k^{1/\alpha} (\log\log T_k)^{\delta-1/\alpha}\right\|-\\
&& \|f\|\left(T_{k+1}^{1/\alpha}(\log\log T_{k+1})^{\delta-1/\alpha}-T_k^{1/\alpha}(\log\log T_k)^{\delta-1/\alpha}\right)-\\
&& M \cdot(1-T_k/T_{k+1})^{1/2}(T_k)^{1/\alpha}\cdot(\log\log T_k)^{\delta-1/\alpha}
\textit{\ \ \ a.s.}
\end{eqnarray*}
\item[$(ii)$]
For\ \  $ T_k= \exp\{k (\log k)^{-3}\}$ the following  holds 
\begin{eqnarray*}
\lim_{k\to \infty}\frac{T_{k+1}^{1/\alpha}(\log\log T_{k+1})^{\delta-1/\alpha}-T_k^{1/\alpha}(\log\log T_k)^{\delta-1/\alpha}}{(T_{k+1}/\log\log T_{k+1})^{1/\alpha}}= 0,
\end{eqnarray*}
\begin{eqnarray*}
\lim_{k\to \infty}\frac{(1-T_k/T_{k+1})^{1/2}(T_k)^{1/\alpha}\cdot(\log\log T_k)^{\delta-1/\alpha}}{(T_{k+1}/\log\log T_{k+1})^{1/\alpha}}= 0,
\end{eqnarray*}
\begin{eqnarray*}
\lim_{k\to \infty}\frac{T_k}{T_{k+1}}= 1.
\end{eqnarray*}\end{itemize}
\end{lemma}
\begin{proof}
To proof $(i)$,
choose  $\tau_k\in [T_k, T_{k+1}]$ such that 
$$ \left\|X_{\alpha}(\tau_k \cdot) - f(\cdot)\tau_k^{1/\alpha} (\log\log \tau_k)^{\delta-1/\alpha}\right\|=\inf_{[T_k, T_{k+1}]} \left\|X_{\alpha}(T\cdot) - f(\cdot)T^{1/\alpha} (\log\log T)^{\delta-1/\alpha}\right\|. $$
Then we need some cumbersome  computations
\begin{eqnarray*}
&&\left\|X_{\alpha}(T_k\cdot) - f(\cdot)T_k^{1/\alpha} (\log\log T_k)^{\delta-1/\alpha}\right\|=
\sup_{s\in[0,T_k]}|X_{\alpha}(s) - f(s/T_k)T_k^{1/\alpha} (\log\log T_k)^{\delta-1/\alpha}|=\\
&&\sup_{s\in[0,T_k/\tau_k]}|X_{\alpha}(s \tau_k) - f(s\tau_k/T_k)T_k^{1/\alpha} (\log\log T_k)^{\delta-1/\alpha}|\leq  \\
&& \sup_{s\in[0,T_k/\tau_k]}|X_{\alpha}(s \tau_k) - f(s)\tau_k^{1/\alpha} (\log\log \tau_k)^{\delta-1/\alpha}|+\\
&& \sup_{s\in[0,T_k/\tau_k]}|f(s\tau_k/T_k)T_k^{1/\alpha} (\log\log T_k)^{\delta-1/\alpha}- f(s)\tau_k^{1/\alpha} (\log\log \tau_k)^{\delta-1/\alpha}|\leq\\
&& \|X_{\alpha}( \tau_k \cdot) - f(\cdot)\tau_k^{\frac{1}{\alpha}} (\log\log \tau_k)^{\delta-1/\alpha}\|+ \sup_{s\in[0,\frac{T_k}{\tau_k}]}|f(s)| \, \left(\tau_k^{\frac{1}{\alpha}} (\log\log \tau_k)^{\delta-1/\alpha}-\right. \\
&& \left. -T_k^{\frac{1}{\alpha}} (\log\log T_k)^{\delta-1/\alpha}\right)+
\sup_{s\in[0,T_k/\tau_k]}|f(s\tau_k/T_k)-f(s)| \,T_k^{1/\alpha} (\log\log T_k)^{\delta-1/\alpha}.
\end{eqnarray*}
Let us show that there exists $0<M<\infty$ such that
$\|f(T_k/\tau_k\cdot)-f(\cdot)\|\leq M \cdot(1-T_k/T_{k+1})^{1/2}$.  
Note that $(1-T_k/T_{k+1})^{1/2}<1$. 

It is known that $\bf{H}$ is dense in $\cal{C}$. 
Thus, for any $k>0$ there exists $f_k \in \bf{H}$ such that $\|f_k(\cdot)-f(\cdot)\|\leq (1-T_k/T_{k+1})^{1/2}$.
We can write 
\begin{eqnarray*}
&& \|f(T_k/\tau_k\cdot)-f(\cdot)\|\leq 2\|f_k(\cdot)-f(\cdot)\| +\|f_k(T_k/\tau_k\cdot)-f_k(\cdot)\|\leq \\
&& (2+\|f'_k\|_{L_2}) \cdot(1-T_k/T_{k+1})^{1/2}.
\end{eqnarray*}
For the last step we used: for any $0 \leq a\leq 1$, any $s\in [0,1]$  we have  $|f_k (a\, s)- f_k(s)|\leq \|f'_k\|_{L_2}(1-a)^{1/2}$, that is easy to prove by Schwarz's inequality.

The rest is obvious.

\medskip

To prove $(ii)$, note
\begin{eqnarray*}
1 \geq \frac{T_k}{T_{k+1}}=\exp \left\{ \frac{k}{(\log k)^3}-\frac{k+1}{(\log (k+1))^3}\right\}\geq \exp \left\{ -\frac{1}{(\log k)^3}\right\},
\end{eqnarray*}
and $\log\log T_k= \log k \,(1+o(1))$. It is left just to make computations.
\end{proof}

Let us show that
\begin{eqnarray*}
\liminf_{k\to\infty}\,\frac{ \left\|X_{\alpha}(T_k\cdot) -f(\cdot)T_k^{1/\alpha} (\log\log T_k)^{\delta-1/\alpha}\right\|}{(T_{k+1}/\log\log T_{k+1})^{1/\alpha}}\geq K_{\alpha}^{1/\alpha}.
\end{eqnarray*}
Take $A>0$. We use the Anderson inequality, self-similatity, and estimate (\ref{Mog}) to obtain
\begin{eqnarray*}
&& \Pr\left\{  \left\|X_{\alpha}(T_k\cdot) -f(\cdot)T_k^{1/\alpha} (\log\log T_k)^{\delta-1/\alpha}\right\|<
A (T_{k+1}/\log\log T_{k+1})^{1/\alpha}\right\}\leq \\
&& \Pr\left\{  \left\|X_{\alpha}(T_k\cdot)\right\|<
A \,(T_{k+1}/\log\log T_{k+1})^{1/\alpha}\right\} \leq \\
&&  \Pr\left\{  \left\|X_{\alpha}(\cdot)\right\|<
A\, \left(\frac{T_{k+1}}{T_k}\right)^{1/\alpha}(\log\log T_{k+1})^{-1/\alpha}\right\} \leq \\
&&\leq \exp\left\{-\frac{K_{\alpha}}{A^{\alpha}}\frac{T_k}{T_{k+1}} \log\log T_{k+1}(1+o(1))\right\}.
\end{eqnarray*}
Using the particular form of $\{T_k\}_{k>0}$, we get
\begin{eqnarray*}
 \Pr\left\{  \left\|X_{\alpha}(T_k\cdot) -f(\cdot)T_k^{1/\alpha} (\log\log T_k)^{\delta-1/\alpha}\right\|<
A\, (T_{k+1}/\log\log T_{k+1})^{1/\alpha}\right\}\leq  \\ \exp\left\{-
\frac{K_{\alpha}}{A^{\alpha}}\log (k /(\log k)^3 )(1+o(1))\right\}.
\end{eqnarray*}
Choose $A=(K_{\alpha}/(1+\epsilon))^{1/\alpha}$, $\epsilon>0$ obtain
\begin{eqnarray*}
 \Pr\left\{  \left\|X_{\alpha}(T_k\cdot) -f(\cdot)T_k^{1/\alpha} (\log\log T_k)^{\delta-1/\alpha}\right\|<
 \left(\frac{K_{\alpha}}{1+\epsilon}\cdot \frac{T_{k+1}}{\log\log T_{k+1}}\right)^{1/\alpha}\right\}\leq  \left(\frac{(\log k)^3}{k}\right)^{(1+\epsilon)(1+o(1))}.
\end{eqnarray*}
 Use the Borel-Cantelli lemma and obtain that for any $f \in \cal{C}$ , any $\epsilon>0$ the following holds
\begin{eqnarray*}
 \liminf_{k\to\infty}\,\frac{ \left\|X_{\alpha}(T_k\cdot) -f(\cdot)T_k^{1/\alpha} (\log\log T_k)^{\delta-1/\alpha}\right\|}{(T_{k+1}/\log\log T_{k+1})^{1/\alpha}} \geq \left(\frac{K_{\alpha}}{1+\epsilon}\right)^{1/\alpha}\ \ a.s.
\end{eqnarray*}
To conclude the proof, tend $\epsilon \to 0$.  

{\it{Addition to Theorem \ref{Chung_deg}: }} The scheme of the proof is the same for ${h(T)}=o((\log\log T)^{-1/\alpha})$ as $T\to \infty$. The difference is just in the first inequality
\begin{eqnarray*}
 \liminf_{T\to\infty}\,  \left\|\frac{X_{\alpha}(T\cdot)}{T^{1/\alpha}h(T)} - f(\cdot)\right\| \geq
 \liminf_{k\to\infty}\,\frac{ \inf_{[T_k, T_{k+1}]} \left\|X_{\alpha}(T\cdot) -f(\cdot)\,T^{1/\alpha}h(T)\right\|}{(T_{k+1}/\log\log T_{k+1})^{1/\alpha}}.
\end{eqnarray*}
Lemma \ref{lemma} could be modified correspondingly.

\bigskip 
  
{\bf{Upper bound for Theorem \ref{LILL} and Theorem \ref{LIL}:}}  Choose $T_k=\exp\{k^{\gamma}\}$, $\gamma>1$. Consider events
\begin{eqnarray*}
D_k(A)=\left\{\left\|\frac{X_{\alpha}(T_k\cdot)}{T_k^{1/\alpha}
(\log\log T_k)^{\delta-1/\alpha}} - f(\cdot)\right\|\leq
A\, \frac{1}{(\log\log T_k)^{\delta}}
\right\},
\end{eqnarray*}
where $0<A<\infty$, $\delta\in (0,1]$. Let us estimate
\begin{eqnarray*}
\Pr\{D_k(A)\}=\Pr\left\{\left\| X_{\alpha}(\cdot) - f(\cdot)(\log\log T_k)^{\delta-1/\alpha} \right\|\leq
A\, \frac{1}{(\log\log T_k)^{1/\alpha}}
\right\}.
\end{eqnarray*}
For $\delta\in (0,1)$, we use the lower bound of Theorem  \ref{small lambda}
\begin{eqnarray*}
&& \Pr\left\{\left\| X_{\alpha}(\cdot) - f(\cdot)(\log\log T_k)^{\delta-1/\alpha} \right\|\leq
A\, \frac{1}{(\log\log T_k)^{1/\alpha}}
\right\}\geq \\ && \exp\left\{- \frac{K_{\alpha}}{A^{\alpha}}\log\log T_k\,(1+o(1))\right\}=
\exp\left\{- \frac{K_{\alpha}}{A^{\alpha}}\gamma \log k\,(1+o(1))\right\}.
\end{eqnarray*}
Putting $A_{\gamma}= (K_{\alpha}\gamma)^{1/\alpha}$, we obtain $\Pr\{D_k(A_{\gamma})\}\geq 1/k$. Thus,
\begin{eqnarray}\label{D}
\sum_{k=1}^{\infty}\Pr\left\{D_k(A_{\gamma})\right\}=\infty.
\end{eqnarray}
For $\delta=1$, we use the lower bound of Theorem  \ref{LIL case} that holds for any 
$\|f'\|<\frac{2}{2-\alpha}\cdot \frac{1}{A^{\alpha-1}}$
and obtain
\begin{eqnarray*}
&&\Pr\left\{D_k(A) \right\} = \Pr\left\{\left\| X_{\alpha}(\cdot) - f(\cdot)(\log\log T_k)^{1-1/\alpha} \right\|\leq
A\, \frac{1}{(\log\log T_k)^{1/\alpha}}
\right\}\geq \\ && \exp\left\{- \frac{C(\alpha)}{A^{\alpha}}\log\log T_k\,(1+o(1))\right\}=
\exp\left\{- \frac{C(\alpha)}{A^{\alpha}}\gamma \log k\,(1+o(1))\right\}
\end{eqnarray*}
Put $A_{\gamma}= (C(\alpha) \gamma)^{1/\alpha}$, and obtain (\ref{D}).

We could not use the Borel-Cantelli lemma directly  because the events
$\{D_k\}$ are dependent. To overcome this difficulty we decompose the process into a sum of independent processes:
\begin{eqnarray}\label{dec}
X_{\alpha}(T_k\cdot)=Y_k(\cdot)+Z_k(\cdot) \textit{\ a.s.},
\end{eqnarray}
where
$$
Y_k(s)=\left\{
\begin{array}{rcr}& X_{\alpha}(T_k s),& s\in [0, \frac{T_{k-1}}{T_k}]\\
& X_{\alpha}(T_{k-1}),& s\in [\frac{T_{k-1}}{T_k}, 1]
\end{array}
\right.
$$
and
$$
Z_k(s)=\left\{
\begin{array}{rcl}& 0,& s\in [0, \frac{T_{k-1}}{T_k}]\\
& X_{\alpha}(T_k s)-X_{\alpha}(T_{k-1}),& s\in
[\frac{T_{k-1}}{T_k}, 1].
\end{array}
\right.
$$
It is easy to see that $Z_1(\cdot), Z_2(\cdot), ... Z_k(\cdot),
Z_{k+1}(\cdot),...$ are independent processes (they are
constructed by using increments of $X_\alpha$ at non-intersecting
intervals).

Let us prove the following
\begin{eqnarray}\label{z}
\sum_{k=1}^{\infty} \Pr\left\{ \left\|\frac{Z_{k}(\cdot)}{T_k^{1/\alpha}
(\log\log T_k)^{\delta-1/\alpha}} -f(\cdot)\right\|\leq
(1+\epsilon)A_{\gamma}\frac{1}{\log\log T_k}\right\}=\infty.
\end{eqnarray}
We use
\begin{eqnarray*}
&&\Pr\left\{ \left\|\frac{Z_{k}(\cdot)}{T_k^{1/\alpha}
(\log\log T_k)^{\delta-1/\alpha}} -f(\cdot)\right\|\leq
(1+\epsilon)A_{\gamma}\frac{1}{(\log\log T_k)^{\delta}}\right\}\geq\\
&&\Pr\left\{ \left\|\frac{X_{\alpha}(T_k\cdot)}{T_k^{1/\alpha}
(\log\log T_k)^{\delta-1/\alpha}} -f(\cdot)\right\|+\frac{\left\|Y_{k}(\cdot)\right\|}{T_k^{1/\alpha}
(\log\log T_k)^{\delta-1/\alpha}}\leq
(1+\epsilon)A_{\gamma}\frac{1}{(\log\log T_k)^{\delta}}\right\}\geq\\
&& \Pr\left[\left\{ \left\|\frac{X_{\alpha}(T_k\cdot)}{T_k^{1/\alpha}
(\log\log T_k)^{{\delta-1/\alpha}}} -f(\cdot)\right\|\leq \frac{A_{\gamma}}{(\log\log T_k )^{\delta}}\right\}\cap
\left\{\frac{\left\|Y_{k}(\cdot)\right\|}{T_k^{1/\alpha}
(\log\log T_k)^{-\frac{1}{\alpha}}}\leq
{\epsilon\,A_{\gamma}}\right\} \right]\geq\\
&& \Pr\{D_k(A_{\gamma})\}- \Pr \left\{\frac{\left\|Y_{k}(\cdot)\right\|}{T_k^{1/\alpha}
(\log\log T_k)^{-1/\alpha}}\geq
\epsilon {A_{\gamma}}\right\}.
\end{eqnarray*}
It is left to prove that the second term could be majorized by a term of convergent series.
Take arbitrary $\epsilon>0$.  Consider the events
\begin{eqnarray*}
C_k(\epsilon A_{\gamma})=\left\{\frac{\left\|Y_{k}(\cdot)\right\|}{T_k^{1/\alpha}(\log\log
T_k)^{-1/\alpha}}>\epsilon\,  A_{\gamma}\right\}.
\end{eqnarray*}
Now we need a large deviation result  (cf. p.238, \cite{Ber_96})
\begin{eqnarray*}
\Pr\{\left\| X_{\alpha}(\cdot)\right\|>x\}=Kx^{-\alpha}(1+o(1))
\textit{ as } x\to\infty,
\end{eqnarray*}
what is true for some $0<K<\infty$. Using
$\left\|Y_{k}(\cdot)\right\|=\left\|X_{\alpha}(T_{k-1}\cdot)\right\|$
a.s. and the self-similarity we write
\begin{eqnarray*}
&&
\Pr\{C_k(\epsilon A_{\gamma})\}=\Pr\left\{\frac{\left\|X_{\alpha}(T_{k-1}\cdot)\right\|}{T_k^{1/\alpha}(\log\log
T_k)^{-1/\alpha}}>\epsilon\,  A_{\gamma}\right\}=\\
&&
\Pr\left\{\frac{\left\|X_{\alpha}(T_{k-1}\cdot)\right\|}{T_{k-1}^{1/\alpha}}>\epsilon\,  A_{\gamma} {\left(\frac{T_k}{T_{k-1}}\right)^{1/\alpha}(\log\log
T_k)^{-1/\alpha}}\right\}=\epsilon\,  A_{\gamma}\,  K \cdot
\frac{T_{k-1}}{T_{k}}\log\log T_k\, (1+o(1)).
\end{eqnarray*}
Now we use
\begin{eqnarray*}
\sum_{k=1}^{\infty}\frac{T_{k-1}}{T_{k}}\log\log T_k=\gamma
\sum_{k=1}^{\infty}\frac{\log k}{\exp\{\gamma
k^{\gamma-1}\}}(1+o(1))<\infty.
\end{eqnarray*}
Thus, for any $\epsilon>0$, any $\gamma>1$ we have
\begin{eqnarray*}
\sum_{k=1}^{\infty}\Pr\{C_k(\epsilon A_{\gamma})\}<\infty.
\end{eqnarray*}
Using Borel-Cantelli lemma we also obtain
\begin{eqnarray*}
\limsup_{k\to \infty}\,(\log\log
T_k)^{1/\alpha}\frac{\left\|Y_{k}(\cdot)\right\|}{T_k^{1/\alpha}}=0.
\end{eqnarray*}
So, this and (\ref{D}) prove (\ref{z})  and the events there are independent. We apply the
Borel-Cantelli lemma and obtain
\begin{eqnarray*}
\liminf_{k\to \infty}\,(\log\log T_k)^{\delta}
\left\|\frac{Z_{k}(\cdot)}{T_k^{1/\alpha} (\log\log
T_k)^{\delta-1/\alpha}} -f(\cdot)\right\|\leq
(1+\epsilon)A_{\gamma}.
\end{eqnarray*}
It is left just to use the elementary relations
\begin{eqnarray*}
&& \liminf_{T\to \infty}\, (\log\log T)^{\delta}
\left\|\frac{X_{\alpha}(T\cdot)}{T^{1/\alpha}
(\log\log T)^{\delta-1/\alpha}} -f(\cdot)\right\|\leq \\
&& \liminf_{k\to \infty}\,(\log\log T_k)^{\delta}
\left\|\frac{X_{\alpha}(T_k\cdot)}{T_k^{1/\alpha}
(\log\log T_k)^{\delta-1/\alpha}} -f(\cdot)\right\|\leq
\\
&& \liminf_{k\to \infty}\,(\log\log
T_k)^{\delta}\left[\left\|\frac{Z_{k}(\cdot)}{T_k^{1/\alpha} (\log\log
T_k)^{\delta-1/\alpha}} -f(\cdot)\right\|+
\frac{\left\|Y_{k}(\cdot)\right\|}{T_k^{1/\alpha}
(\log\log T_k)^{\delta-1/\alpha}}\right]\leq \\
&&\liminf_{k\to \infty}\,(\log\log T_k)^{\delta}
\left\|\frac{Z_{k}(\cdot)}{T_k^{1/\alpha} (\log\log
T_k)^{\delta-1/\alpha}} -f(\cdot)\right\|+ \limsup_{k\to
\infty}\,(\log\log
T_k)^{1/\alpha}\frac{\left\|Y_{k}(\cdot)\right\|}{T_k^{1/\alpha}}.
\end{eqnarray*}
Tending $\gamma\to 1$ and $\epsilon \to 0$, we obtain the upper bound.
\end{proof}

{\bf{Open questions:}}
\begin{enumerate}
 \item Wide field of action is to find a.s. limit sets  in the case of scaling  functions 
 \begin{eqnarray*}
 (\log\log T)^{1-1/\alpha}<h(T)\leq  (\log T)^{1/\alpha}(\varphi(\log\log T))^{1/\alpha}, 
 \end{eqnarray*}
 where $\varphi$ is as in (\ref{up_norm}). From the proof we see that a positive result (the a.s. limit set  is wider than $\{\mathbf{0}\}$) requires a good lower bound of $\Pr\{\|X_{\alpha}-\lambda f\|< r\}$ under $\lambda r^{\alpha-1}\to \infty$, and a negative result (the a.s. limit set  is $\{\mathbf{0}\}$) would require a good upper bound of the same probability. 
 \item It is interesting to study the functional LIL in the Skorokhod topology. We already mention that the a.s. cluster set will contain the a.s. cluster set under the uniform convergence. It is also possible that the set of admissible scaling functions is wider.
  \end{enumerate}
\section*{Acknowledgment}
The author is grateful to the department of Financial and
Actuarial Mathematics of the Vienna Technical University, in particular, to
Reinhold Kainhofer, Josef Teichmann  and Friedrich Hubalek for the
discussions and the encouragement and also to Mikhail Lifshits for showing a parallel to Baldi-Royonette results.

\bibliographystyle{alpha}
\bibliography{mybibliography}

\end{document}